\documentclass[12pt]{article}
\usepackage{amsmath}
\usepackage{amscd}
\usepackage{amssymb}
\usepackage{array}

\begin{document}

\newcommand{\R}{\mathbb{R}}
\newcommand{\C}{\mathbb{C}}
\newcommand{\Sc}{\mathcal{S}}
\newcommand{\Pol}{\rm Pol}
\newcommand{\lra}{\longrightarrow}
\newcommand{\s}{\sigma}
\newcommand{\fh}{\mathfrak{h}}
\newcommand{\fg}{\mathfrak{g}}

\newcommand{\B}{\mathcal{B}}
\newcommand{\G}{\mathcal{G}}
\newcommand{\A}{\mathcal{A}}
\newcommand{\F}{\mathcal{F}}
\newcommand{\I}{\mathcal{I}}
\newcommand{\J}{\mathcal{J}}
\newcommand{\Ha}{\mathcal{H}}
\newcommand{\Oc}{\mathcal{O}}
\newcommand{\La}{\mathcal{L}}

\newcommand{\m}{\mathbf{m}}
\newcommand{\Hb}{\mathbf{H}}
\newcommand{\Gb}{\mathbf{G}}
\newcommand{\Sym}{{\rm {Sym}}}
\newcommand{\N}{\mathcal{N}_{\Theta}}
\newcommand{\U}{\mathcal{U}}
\newcommand{\de}{\delta}
\newcommand{\Z}{\mathbb{Z}}
\newcommand{\sym}{{\rm Sym}}
\newcommand{\Gc}{G_{\C}}
\newcommand{\gc}{\g_{\C}}
\newcommand{\la}{\lambda}
\newcommand{\Ss}{\mathcal{S}}

\newcommand{\h}{\hbar}
\vskip 2cm
  \centerline{\LARGE \bf A Review on Deformation Quantization of}
  \bigskip
\centerline{\LARGE\bf   Coadjoint Orbits of Semisimple Lie Groups}

\vskip 2.5cm

\centerline{R. Fioresi$^\ast$\footnote{Investigation supported by
the University of Bologna, funds for selected research topics.}
 and M. A. Lled\'o$^\dagger$}

\bigskip

\centerline{\it $^\ast$Dipartimento di Matematica, Universit\`a di
Bologna }
 \centerline{\it Piazza di Porta S. Donato, 5.}
 \centerline{\it40126 Bologna. Italy.}
\centerline{{\footnotesize e-mail: fioresi@dm.UniBo.it}}

\bigskip

\centerline{\it $^\dagger$ Dipartimento di Fisica, Politecnico di
Torino,} \centerline{\it Corso Duca degli Abruzzi 24, I-10129
Torino, Italy, and} \centerline{\it INFN, Sezione di Torino,
Italy.} \centerline{{\footnotesize e-mail:
lledo@athena.polito.it}}
 \vskip 2cm

\begin{abstract}

In this paper we make a review of the results obtained in previous
works by the authors on deformation quantization of coadjoint
orbits of semisimple Lie groups.
We motivate the problem with a new point of view of the well
known Moyal-Weyl deformation quantization.
We consider only semisimple
orbits. Algebraic and differential deformations are compared.

\end{abstract}

\vfill\eject

\section{Motivation}

The difference between classical and quantum mechanics is the
presence of a non commutative algebra substituting the commutative
algebra of classical observables. This non commutative algebra has
a representation in some Hilbert space. The approach of
deformation quantization consists in obtaining the quantum algebra
as a deformation of the classical one. The parameter of
deformation, which  measures the non commutativity of the algebra,
is the  Planck constant $\hbar$. This is in agreement with the
Bohr correspondence principle which states that when $\hbar \lra
0$ the behaviour of the quantum system becomes closer to the
classical one.

To quantize a physical system defined on $\R^2$ one associates to
the coordinate functions $p$ and $q$ on $\R^2$ the self-adjoint
operators on $P$ and $Q$ on $\Ha=L^2(\R)$:
$$ P\psi(q)=-i\hbar
\partial_q\psi(q), \qquad Q\psi(q)=q\psi(q), \qquad \psi \in L^2(\R) $$
satisfying the commutation rule $[Q,P]=-{i\hbar}I$. For  generic
polynomials in $p$ and $q$ one should select an ordering rule.
Using the symmetric (or Weyl) ordering rule the operator
associated to the polynomial $q^np^m$ is self-adjoint and given by
\cite{we}: $$ w(x_{i_1} \cdots x_{i_{n+m}})=\frac{1}{(n+m)!}
\sum_{\sigma \in S_{n+m}} X_{i_{\s(1)}} \dots X_{i_{\s(n+m)}} $$
with $i_j=1,2$ and $x_1=q$, $x_2=p$, $X_1=Q$, $X_2=P$. This
defines a non commutative algebra structure on the polynomial
algebra on $\R^2$, $\Pol(\R^2)$,
\begin{equation}
a \star b=w^{-1}(w(a) \cdot w(b)), \qquad a,b \in \Pol(\R^2)
\label{weylprod}
\end{equation}
whose associated Lie bracket was first written by Moyal \cite{mo}.

This product can be expressed as an infinite power series in
$\hbar$: $$ f \star g=fg+B_1(f,g)\hbar+B_2(f,g)\hbar^2+ \dots
\label{mw} $$ with the $B_i$  bioperators that coincide with
bidifferential operators  acting on polynomials. This means that
one can extend $B_i$ to act on $C^{\infty}(\R^2)$. The degree of
the bidifferential operators is $(i,i)$ (we will see the explicit
form later), which means that the product of two polynomials gives
another polynomial. When extending to $C^\infty(\R^2)$ the product
of two functions may be an infinite, non convergent power series
in $\h$. The product is well defined only in the space of formal
power series in $\hbar$ with coefficients in $C^\infty(\R^2)$.
This is usually denoted by $C^\infty(\R^2)[[\hbar]]$. This was the
approach to differential deformations taken in Ref.\cite{bffls}.

\medskip
It is our intention to give another point of view on the product
defined above by $w$ on $\Pol(\R^2)$ and its generalization to the
$\C^{\infty}$ functions.

\medskip

Let $\Hb$ be the three parameter group $\Hb =\R^3$ with
multiplication $$(a_1,b_1,c_1)\cdot(a_2,b_2,c_2)=
(a_1+a_2,b_1+b_2,c_1+c_2+a_1b_2).$$ It is the Heisenberg group
with Lie algebra $\fh=\R^3={\rm span}\{Q,P,E'=-iE\}$ and
commutation rules $$[Q,P]=-iE\qquad \hbox{(the rest trivial)}.$$

We consider the dual space $\fh^*$ with coordinates $(q,p,e'=-ie)$
in the basis dual to $\{Q,P,E'=-iE\}$. The algebra of real polynomials
on $\fh^*$ is
$${\rm Pol}(\fh^*)={\rm span}_{\R}\{q^mp^n{e'}^r,\quad
m,n,r=0,1,2,\dots\}.$$

There is a Poisson structure on $C^\infty(\fh^*)$
\begin{equation}\{f_1, f_2\}=e'(\frac{\partial f_1}{\partial
q}\frac{\partial f_2}{\partial p}-\frac{\partial f_2}{\partial
q}\frac{\partial f_1}{\partial p})\qquad f_i\in
C^\infty(\fh^*),\label{pbh}\end{equation} under which the
polynomials ${\rm Pol}(\fh^*)\subset C^\infty(\fh^*)$ are closed.
This is the Kirillov Poisson structure for $\fh$.

The Poisson structure is tangent to the planes  $e'={\rm
constant}$, where it restricts  as a symplectic structure (maximal
rank) provided $e'\neq 0$. On the plane $e' =0$ the Poisson structure
is identically 0. The planes $e'={\rm constant}\neq 0$ and  the
points of the form $(q,p,0)$ are the leaves of the symplectic
foliation of the non regular Poisson structure (\ref{pbh}).

\medskip

The adjoint representation of $g=(a,b,c)\in \Hb$ in the ordered
basis $\{Q,P,E'\}$ is $${\rm
Ad}_g=\begin{pmatrix}1&0&0\\0&1&0\\-b&a&1\end{pmatrix}.$$ On
$\fh^*$  we have  the coadjoint action defined by given by $${\rm
Ad}^*_g\xi(X)=\xi({\rm Ad}_{g^{-1}}X) \qquad \xi\in \fh^*, \;
X\in\fh,\; g\in \Hb. $$ In the dual basis, $${\rm
Ad}^*_g=\begin{pmatrix}1&0&b\\0&1&-a\\0&0&1\end{pmatrix}.$$ The
orbits of the coadjoint action are the planes $e'={\rm
constant}\neq 0$ and the single points $(q,p,0)$. They coincide
with the leaves of the symplectic foliation of (\ref{pbh}). (This
is also the case for a general Lie algebra).

So we have that the coadjoint orbits are symplectic manifolds with
an action of the group $\Hb$. Since the stability group (for
$e'\neq 0$ is $\R\approx \{(0,0,c)\}$,  they are in fact coset
spaces $\Hb/\R$.

\smallskip

Let $\fh_\h$ be the Heisenberg algebra with bracket multiplied by
$\h \in \R$ (Later on we will take it as formal parameter). We
consider now the enveloping algebra of $\fh_\h$, $U(\fh_\h)$,
which is the tensor algebra $T(\fh_\h)$ modulo the ideal generated
by the relations of the Lie bracket, that is $$Q\otimes P-P\otimes
Q =\h E',\quad Q\otimes E'-E'\otimes Q=0,\quad P\otimes E'-
E'\otimes Q =0.$$

An ordering rule in $U(\fh_\h)$ is a linear bijection ${\rm
Pol}(\fh_\h^*)\rightarrow U(\fh_\h)$. We can take for example the
symmetric (or Weyl) ordering rule
\begin{equation}
W(x_{i_1}\cdots x_{i_p})={\frac{1}{p!}}\sum_{s \in
S_p}X_{i_{s(1)}}\otimes\cdots\otimes X_{i_{s(p)}},\label{weyl}
\end{equation}
where $x_1=q$, $x_2=p$, $x_3=e'$ and $X_1=Q$,
$X_2=P$, $X_3=E'$ and $i_j = 1,2,3$.

We can define an associative, non commutative product on ${\rm
Pol}(\fh_\h^*)$ as \begin{equation}f_1\star
f_2=W^{-1}(W(f_1)W(f_2)).\label{star}\end{equation} Explicitly,

\begin{equation}f_1\star f_2(q,p,e')=\sum_{k=0}^{\infty}
{\frac{1}{k!}}({\frac{\h}{2}e'})^kP^k(f_1,f_2)=
\exp({-{\frac{i\h}{2}e'P}})(f_1,f_2),\label{sph}\end{equation}
$$P^k(f_1,f_2)=P^{i_1j_1}\cdots P^{i_kj_k}\frac{\partial
f_1}{\partial x^{i_1}\cdots \partial x^{i_k}}\cdot\frac{\partial
f_2} {\partial x^{j_1}\cdots
\partial x^{j_k}},$$ $$P=\begin{pmatrix} 0&1\cr
-1&0\end{pmatrix},\qquad i_l,j_l=1,2, \qquad l=1 \dots k.$$

The product defined in (\ref{sph}) is of the general form $$f\star
g=\sum_{i=0}^\infty C_i(f,g)h^i$$ with  $C_i$  bidifferential
operators and satisfying the properties
\begin{eqnarray}&&1. \lim_{\h\mapsto 0}f_1\star f_2=f_1
f_2,\nonumber\\&&2.\lim_{\h\mapsto 0}\frac{1}{2\h}(f_1\star
f_2-f_2\star f_1)=\{f_1,f_2\}.\label{prop}\end{eqnarray}

Since in (\ref{sph}) only derivatives $\frac{\partial}{\partial
q}$ and $\frac{\partial}{\partial p}$ appear, the bidifferential
operators are tangent to the orbits $e={\rm constant}\neq 0$. One
can in fact restrict it to the orbit (say, $e=1$) without
ambiguity and obtain a non commutative associative product on
$\Pol(\R^2)$. It is the Moyal-Weyl product on $\Pol(\R^2)$, that
we have described in (\ref{weylprod}). The differential property
allows us to extend the product to $C^\infty(\R^2)$.

In general, given a real manifold $M$, a product defined on
$C^\infty(M)[[h]]$ satisfying the properties ({\ref{prop}) is a
called a {\it star product} on $M$ (or on $C^{\infty}(M)$). If the
$C_i$'s are bidifferential we say that the star product is {\it
differential}. If $M$ is an affine algebraic variety, and the star
product is defined on polynomials we say that it is an {\it
algebraic star product} on $M$ (or on $\C[M]$).

The Moyal-Weyl star product on $\R^2$ is differential, when
restricted to $\Pol(\R^2)$ is algebraic  and additionally it
converges on polynomials for all real values of $\h$.

\medskip

We want now to generalize this construction to the coadjoint
orbits of a semisimple Lie algebra. The simplest non trivial
example of a coadjoint orbit is the sphere $S^2\approx {\rm
SU}(2)/{\rm U}(1)$. We can suitably generalize the algebraic
construction using the enveloping algebra of
$\mathfrak{su}(2)=$Lie(SU(2)). This was done in Refs.\cite{fl,ll},
and we will make a review of it in Section \ref{alge} (another
construction for certain types of algebraic manifolds can be seen
in Ref.\cite{ko2}). The other approach is to generalize the
differential construction. This was done for arbitrary symplectic
(in fact, for regular Poisson) manifolds in Refs.\cite{omy,dl,fe}.
For arbitrary Poisson manifolds it was done first in
Ref.\cite{ko}, and with different approaches in
Refs.\cite{ta,cft}.

How can one write a differential star product on an arbitrary,
symplectic manifold $M$? Let $U$ be an open set in $M$ where
Darboux (canonical) coordinates exist, so $U\approx \R^{2n}$. We
can write a star product on $U$ using the Darboux coordinates and
(\ref{sph}). Formula (\ref{mw}) is not invariant under symplectic
transformations (unless they are linear symplectic
transformations), but two star products given for different open
sets $U,V$ (denoted by $\star_U$ and $\star_V$) are isomorphic in
the intersection $U\cap V$.  It is a general principle that one
can ``glue" the star products defined in every open set of a
covering of $M$ into a globally defined star product on $M$,
$\star_M$ if the two following conditions are satisfied:

\smallskip

\noindent 1. In the intersection $U\cap V$ there is an isomorphism
of algebras \begin{eqnarray*}&\phi_{UV}:C^\infty(U)\longrightarrow
C^\infty(V)\\ &\phi_{UV}(f_1\star_U f_2)= \phi_{UV}(f_1)\star_V
\phi_{UV}(f_2), \qquad f_1,  f_2\in C^\infty(U).\end{eqnarray*}

\smallskip

\noindent 2. In the intersection of three open sets $U\cap W\cap
V$ the cocycle condition $$\phi_{UV}=\phi_{UW}\circ \phi_{WV}$$ is
satisfied.

\smallskip

The first condition is satisfied in the symplectic case. The
second condition can also be satisfied since it has been shown
that a global star product exists \cite{omy,dl,fe}. In Appendix 1
we will write an explicit form of the gluing in terms of a
partition of unity.

It is worth mentioning that in some non trivial
cases global Darboux coordinates
exist for the whole manifold.
(This is happening for example in the
case of coadjoint orbits of $SO(2) \times GL_+(2)$, Ref. \cite{jc}).
In this case a star
product is defined on $C^{\infty}(M)[[h]]$ using the
Moyal-Weyl deformation with Darboux coordinates.

\medskip

In the case of the Heisenberg group the star product on $\fh^*$ is
tangent to the orbits so that the restriction to  the orbit is
well defined. For a more general group this is not guaranteed
\cite{cgr}, and we will see that in the semisimple case one cannot
have a star product that is simultaneously algebraic and
differential, that is a star product defined on the polynomials on
the orbit via bidifferential operators. Nevertheless, we can ask
if there is any relation between the star products constructed
using the algebraic method and the differential one. This question
was addressed in Refs.\cite{fll,fl2}. We will review it in Section
\ref{comparison}.

\section{Algebraic construction of a star product \label{alge}}

Let $\Gb$  a real semisimple Lie group and $\fg={\rm Lie}(\Gb)$.  Let
$\{X_1,\dots X_n\}$ be a basis of $\fg$ and $(\xi_1,\dots \xi_n)$
the coordinates on $\fg^*$ with respect to the dual basis. On
$\fg^*$ we have the Kirillov Poisson bracket
\begin{equation}\{f_1,f_2\}(\xi)=\langle[df_1,df_2],\xi\rangle=\sum_{ijk}c_{ij}^k\xi_k\frac{\partial
f_1}{\partial\xi_i}\frac{\partial f_2}{\partial\xi_j}, \qquad f_1,
f_2 \in C^{\infty}(\fg^*) \label{pbg}\end{equation} where
$c_{ij}^k$ are the structure constants of $\fg$. Since $\fg$ is
semisimple, one can identify  $\fg\approx \fg^*$ by means of the
invariant Cartan-Killing form. From now on we will work in the
complex field, so we will take the complexification of an algebra
when needed. The deformed algebras that we will obtain have a
conjugation and the space of fixed points of a conjugation is a
real algebra. We will denote ${\rm Pol}(\fg^*)_\C= \C[\fg^*]$. If
we multiply the structure constants of the algebra $\fg$ by a
parameter $h$, when we set $h$ equal to a real number we obtain an
algebra isomorphic to $\fg$. Let us now treat $h$ as a formal
parameter. Let $\fg_h$ denote the Lie algebra over $\C[[h]]$
obtained from $\fg$ by multiplying the structure constants by the
formal parameter $h$. Let $U_h(\fg)$ (or for brevity $U_h$) denote
its universal enveloping algebra.

 The Weyl map $W:\C[\fg^*][[h]]\rightarrow U_h$ is defined as
 in (\ref{weyl}) and it is a linear bijection. It defines a star
 product on $\fg^*$ via formula (\ref{star}). It is easy to see that this star product is differential.
  There are other
 linear bijections between $\C[\fg^*]$ and $U(\fg)$ that define other
 (equivalent) star products. We
 will make use of this freedom later. But $W$ has the following
 property  (see
Ref. \cite{va} pg. 183):  let $A$ be  an automorphism (derivation)
of $\fg$. It extends to an automorphism (derivation) of $U_h$
denoted by $\tilde A$. It also extends to an automorphism
(derivation) $\bar A$ of  $S(\fg)_\C[[h]]\simeq \C[\fg^*][[h]]$
($S(\fg)_\C$ denotes the complexification of the symmetric tensors
over $\fg$). Then
\begin{equation}W\circ\bar A=\tilde A\circ
W.\label{intertwin}\end{equation} This intertwining property will
be useful for the construction of the deformation on the orbits.

\medskip

The leaves of the symplectic foliation of (\ref{pbg}) coincide
with the coadjoint orbits of $\Gb$. We consider the complexified
group $\Gb_\C$, and a semisimple element $X\in \fg$ (that is, X
has no nilpotent part in the Jordan decomposition). The coadjoint
orbit of this element is an affine algebraic variety  defined over
$\R$. The real form of this orbit  is a union of real orbits. If
$\Gb$ is compact  then the real form of the complex orbit is the
real orbit itself.
We will construct a deformation quantization of the real
form of the complex orbit.

In the symplectic  foliation of $\fg^*$ there are leaves of
different dimensions. When the dimension is maximal we say that
the orbit (or leaf) is regular. In this case, the polynomials that
define the affine variety are invariant polynomials. In fact, they
are of the form $$p_i=c_i^0, \quad i=1,\dots m={\rm
rank}(\fg),\quad c_i^0\in \R,$$ where $p_i$ are homogeneous
polynomials that generate the whole algebra of invariant
polynomials (Chevalley theorem) $$ {\rm Inv}(\fg)_\C=\C[p_1,\dots
,p_m]. $$ $c_i^0$ are generic constants and $\{dp_i\}$ are
independent on the regular orbit \cite{va2}. The ideal of the
orbit is $\I_0=(p_i-c_i^0, i=1 \dots m)$, and the ring of
polynomials on the orbit is $$ \C[\Theta]=\C[\fg^*]/\I_0. $$

If the orbit is not regular, then one can show \cite{ll} that the
ideal of the orbit is generated by a set of polynomials
$\{r_\alpha,\; \alpha=1,\dots ,l\}$ satisfying
\begin{equation}r_\alpha(g^{-1}X)=\sum_{\beta=1}^lT(g)_{\alpha\beta}r_\beta(X),
\label{rep}\end{equation}
where $T$ is a finite dimensional
representation of $\Gb$. This means that if not the
polynomials, the set itself is invariant under the action of the
group

\medskip

It is not difficult to see that the star product defined on
$\fg^*$ via the Weyl map $W$ is not tangent to the orbits, so the
restriction is not well defined. But there are other isomorphisms
$\psi:\C[\fg^*][[h]]\rightarrow U_h$ and star products
$$f_1\star_\psi f_2=\psi^{-1}(\psi(f_1)\psi(f_2)).$$ We are going
to construct an isomorphism  $\psi$ such that the star product
$\star_\psi$ is well defined on the orbits.

Let $R_\alpha=W(r_\alpha)$. Notice that (\ref{rep}) includes as a
particular case the regular orbits, since one can take
$\{r_\alpha\}=\{p_i-c_i^0\}$. We consider now the two sided ideal
in $U_h$ generated by $R_\alpha$, $\I_h= (R_\alpha)$. It is easy
to show using property (\ref{intertwin}) that the left and right
ideals are equal, and then equal to $\I_h$, $$\I_h^{\rm
left}=\I_h^{\rm right}=\I_h.$$ We look for an isomorphism $\psi$
such that $\psi(\I_0)=\psi(\I_h)$ so the diagram
\begin{equation}
\begin{CD}
\C[\fg^*][[h]]@>\psi>>U_{h}\\ @VV{\pi}V @VV{\pi_h}V\\
\C[\fg^*][[h]]/\I_0@>\tilde\psi>>U_{h}/\I_{h}
\end{CD}
\label{cd}
\end{equation}
commutes and the induced map $\tilde\psi$ is an isomorphism. Then
the product defined by:
$$\tilde f_1\star_{\tilde\psi}\tilde
f_2=\tilde\psi^{-1}(\tilde\psi(f_1)\tilde\psi(f_2))$$ is an
algebraic  star product on the orbit.

In \cite{fl} the case of regular orbits was solved and in
\cite{ll} the regularity condition was removed. Here we show
explicitly the case of SU(2), although the construction was done
in complete generality.

\subsection{Star products on the coadjoint orbits of SU(2)}

We consider the Lie algebra $\fg=\mathfrak{su}(2)={\rm
span}_{\R}\{X,Y,Z\}$ with commutation rules $$[X,Y]=Z,\quad
\hbox{and cyclic permutations}.$$ Let $\{x,y,z\}$ be coordinates
on $\fg^*$ with the dual basis. The coadjoint orbits of SU(2) are
spheres $S^2$ centered in the origin. As affine algebraic
varieties they are defined by the polynomial constrain $$p\equiv
x^2+y^2+z^2=r^2, \qquad r\in \R.$$ The origin itself is the only
non regular orbit.

We now go to the complexification. Let $\I_0=(p-r^2)
\subset \C[x,y,z]$.
Consider the following basis of $\C[x,y,z]$,
\begin{eqnarray*}&&B=B_0\cup B_1\\&& B_0=\{ x^my^nz^q(p-r^2), \;\;
 m,n,q=0,1,2,\dots\}\\&&B_1=\{x^my^nz^\nu,\;\;
\nu=0,1, \; m,n=0,1,2,\dots\}.\end{eqnarray*} $B_0$ is a basis of $\I_0$ and
the set of equivalence classes in $\C[\fg^*]/\I_0$ of the elements
of $B_1$ is a basis of $\C[\Theta]$.  Consider the map
\begin{eqnarray*}&&\psi(x^my^nz^q(p-r^2))=X^mY^nZ^q(P-r^2)\\&&
\psi(x^my^nz^\nu)=X^mY^nZ^\nu.\end{eqnarray*} where
$P=X^2+Y^2+Z^2$ is the Casimir operator of $\mathfrak{su}(2)$.
Then $\psi(\I_0)=\I_h$ with
$\I_h=(P-r^2) \subset U_h$ and the induced map $\tilde \psi$ is an
isomorphism, so an algebraic star product on $S^2$ is defined. In
\cite{fll} it was proven that it is not differential. Notice also
that the star product is tangential only to the orbit with radius
$r$, and not to the orbit with radius $r+\delta r$.

In the general case it is tricky to show that the images under
$\psi$ of $B_1$ are linearly independent in $U_h/\I_h$. In
Ref.\cite{fl} this was done for regular orbits and in
Ref.\cite{ll} we were able to remove the regularity hypothesis.

\medskip

There are other isomorphisms that one can use. In Ref.\cite{cg} it
was explicitly constructed  a star product on the regular orbits
by finding another such isomorphism. The construction is based in
the decomposition \cite{kos} $${\rm Pol}(\fg^*)={\rm
Inv}(\fg^*)\otimes {\rm Harm}(\fg^*).$$ ${\rm Inv}(\fg^*)$ are the
polynomials on $\fg^*$ that are invariant under the action of the
group and ${\rm Harm}(\fg^*)$ are the harmonic polynomials, which
can be identified with the polynomials on the orbit, ${\rm
Pol}(\Theta)$. Then, a  basis of ${\rm Pol}(\fg^*)$ is
$$\{p^rf,\quad f\in {\rm Harm}(\fg^*), \; r=0,1,2,\dots\}.$$ The
isomorphism is given in terms of this basis as
\begin{equation}\psi(p^rf)=P^rW(f).\label{harm}\end{equation}
 The corresponding star product is not
differential \cite{cg}, but it has the property that it is tangent
to all the orbits in a neighborhood of the orbit $p=r^2$.

\medskip

We want to make the observation that in the case of regular
orbits, being $p$ an  invariant polynomial, we can slightly
generalize  the construction and choose $\psi(p-r^2)=P-c(h)$ where
$c(0)=r^2$. This allows to chose $c(h)=r(r+h)$. For each
irreducible representation of  $\mathfrak{su}(2)$ and its
enveloping algebra we can find special values of $r$ in such a way
that the representation descends to   the quotient $U_h/\I_h$. The
image of this algebra under the irreducible representation is a
finite dimensional algebra that was studied in Ref.\cite{gb}. It
has been recently called a ``fuzzy" sphere in the physics
literature. It is also the algebra of geometric quantization
\cite{fl}.

In this way we can construct a family of algebraic deformations.
It turns out that they are not all isomorphic  \cite{fll}.

\section{Comparison between the algebraic and the differential
methods for regular orbits\label{comparison}}

We have seen in the case of SU(2) how to construct a family of
deformations of the regular orbits. The construction can be
generalized to all the regular orbits of a compact semisimple Lie
group. (For the non compact case and semisimple orbits, we recall
that the algebraic variety is indeed a union of orbits). The
deformations are of the form $U_h/\I_h$ where $\I_h=(P_i-c_i(h))$
with $P_i=W(p_i)$, the Casimir operators of $\fg$. Different
choices for $c_i(h)$ (but always with $c_i(0)=c_i^0$) may give non
isomorphic algebras. The question is if any of these algebras is
somehow equivalent to a differential one. The first thing we note
is that we may perhaps get an injective homomorphism of
the star product algebra of polynomials into the $C^\infty$
functions with a differential star product algebra, so
``equivalence" in this context will mean to have such embedding.

If $\fg$ is a semisimple Lie algebra, it is not possible to have a
differential star product that is tangential to all the orbits of
$\fg^*$ \cite{cgr}. Instead, we can consider a regularly foliated
neighborhood of the orbit, where a differential, tangential star
product always exists \cite{ma}.

We will consider three different star products:

\subparagraph{Star product $\star_S$ on $\fg^*$.} It is the one
induced by the Weyl map. It is algebraic, differential, defined on
all $\fg^*$ and {\sl not tangential}.

\subparagraph{Star products $\star_P$ and $\star_{P\Theta}$.} It
is the star product on $\fg^*$ defined by
$$\psi((p_1-c_1^0)^{q_1}\cdots (p_m-c_m^0)^{q_m} f)=
(P_1-c_1(h))^{q_1}\cdots (p_m-c_m(h))^{q_m} W(f),$$ where $c_i(h)$
is still to be determined and $f$ is an harmonic polynomial (see
(\ref{harm})). It is algebraic, {\sl not differential}, defined on
all $\fg^*$ and tangential. We will denote $\star_{P\Theta}$ the
restriction to the orbit.

\subparagraph{Star products $\star_T$ and $\star_{T\Theta}$.} We
consider $\N$, a regularly foliated neighborhood of the orbit, and
an atlas given by a good covering with Darboux coordinates. On
each open set, we construct the star product given by Kontsevich's
local formula \cite{ko} in the Darboux coordinates. From
Kontsevich's theorem, it follows that it is equivalent to the
restriction of $\star_S$ to that open set. The star products
obtained with the Darboux coordinates can be glued into a global
star product on $\N$ (see Appendix 1).

It is differential, tangential and {\it not algebraic}. Also, it
is  equivalent to $\star_S$ on $\N$. Its restriction to $\Theta$
is denoted by $\star_{T\Theta}$.

\bigskip

\bigskip

We have the following maps:
\begin{equation*}
\begin{CD}
(\C[\fg^*][[h]],\star_P)@>\eta>\approx>
(\C[\fg^*][[h]],\star_S)\subset(C^\infty(\N)[[h]],\star_S)\\
(C^\infty(\N)[[h]],\star_S)@>\rho>{\rm
injective}>(C^\infty(\N)[[h]],\star_T).
\end{CD}
\end{equation*}
So we have that $\rho\circ\eta$ is an injective homomorphism.
Since both, $\star_P$ and $\star_T$ are tangential, we may ask if
$\rho\circ\eta$ descends as an injective homomorphism to the
algebras on the orbits. In \cite{fl2} it was shown that $c_i(h)$
can be chosen in such a way that $\rho\circ\eta(\I_0)=\I_0$ so we
have that the polynomial algebra
$(\C[\Theta][[h]],\star_{P\Theta})$ inject homomorphically into
$(C^{\infty}(\Theta)[[h]],\star_{T\Theta})$.

\section*{Appendix 1. Gluing of star products \label{glue}}

Let $M$ be a Poisson manifold and fix an open cover $\U=\{U_r\}_{r
\in J}$ where $J$ is some set of indices. Assume that in each
$U_r$ there is a differential star product $$
\star_r:C^\infty(U_r)[[h]]\otimes
C^\infty(U_r)[[h]]\longrightarrow C^\infty(U_r)[[h]] $$ and
isomorphisms
\begin{eqnarray*}
&&T_{sr}:C^\infty(U_{rs})\longrightarrow C^\infty(U_{sr}),\qquad
U_{r_1\dots r_k}=U_{r_1}\cap \cdots \cap U_{r_k} \nonumber\\
&&T_{sr}(f)\star_s T_{sr}(g)=T_{sr}(f\star_rg)\end{eqnarray*} such
that the following conditions are satisfied
\begin{eqnarray}
&1.& T_{rs}=T_{sr}^{-1} \qquad\qquad\; {\rm{on}} \;
U_{sr},\nonumber
\\ &2.& T_{ts}\circ T_{sr}=T_{tr}, \qquad\, {\rm{on}}\;
U_{rst}.\label{cocyco}\end{eqnarray} Then there exists a global
star product on $M$ isomorphic to local star product on each
$U_r$.

\medskip

Let $\phi_i:U_i\rightarrow \R$ be a partition of unity of $M$
subordinate to the covering $\U$ and $f_r\in C^\infty(U_r)$ such
that $f_r=T_{rs}f_s.$ Let $f\in C^\infty(M)$ be the function
$$f=\sum_{r\in J}\phi_rf_r. $$ On $U_r$, $f$ becomes
$$f=(\phi_r{\rm{Id}}+\sum_s\phi_sT_{sr})f_r=A_rf_r.$$ We define a
star product on $U_r$  as \begin{equation}f\star
g=A_r(A_r^{-1}(f)\star_rA_r^{-1}(g)).\label{mgsp}\end{equation} it
is equivalent to $\star_r$.  Using conditions (\ref{cocyco}) one
has $$A_rT_{rt}=A_t.$$ Then, the star products
 (\ref{mgsp}) on each $U_r$ coincide in the intersections, so they
define a unique star product on $M$. The restriction  of this star
product to $U_r$ is equivalent to $\star_r$.  Also, using
different partitions of unity one obtains equivalent star
products.

\end{document}